\documentclass[12pt]{amsart}
\usepackage{amssymb}
\usepackage{amscd}
\usepackage{pstricks}
\usepackage{pst-node}
\newtheorem{thm}{Theorem}[section]
\newtheorem{cor}[thm]{Corollary}

\newtheorem{prop}[thm]{Proposition}

\theoremstyle{definition}

\newtheorem{rem}[thm]{Remark}

\theoremstyle{remark}

\numberwithin{equation}{section}

\newcommand{\trans}[1]{{}^t\kern-.2em{#1}}
\newcommand{\ytrans}[1]{{}^t\kern-.11em{#1}}
\newcommand{\Trans}[1]{{}^T\kern-.2em{#1}}
\newcommand{\lsup}[2]{{}^{#1}\kern-.1em{#2}}

\newcommand{\Id}{\operatorname{Id}}

\newcommand{\Ker}{\operatorname{Ker}}

\newcommand{\Mas}{\mathbf{Mas}}
\newcommand{\M}{\mathbf{\M}}
\newcommand{\Sf}{\mathbf{Sf}}

\renewcommand{\tilde}[1]{\widetilde{#1}}

\DeclareFixedFont{\bgn}{OT1}{cmr}{m}{n}{20.74}
\DeclareFixedFont{\bgi}{OT1}{cmr}{m}{it}{20.74}

\newcommand{\bigzerou}{\smash{\lower1.7ex\hbox{\bgi O}}}

\makeatletter
\def\eqnarray{%
   \stepcounter{equation}%
   \def\@currentlabel{\p@equation\theequation}%
   \global\@eqnswtrue
   \m@th
   \global\@eqcnt\z@
   \tabskip\@centering
   \let\\\@eqncr
   $$\everycr{}\halign to\displaywidth\bgroup
       \hskip\@centering$\displaystyle\tabskip\z@skip{##}$\@eqnsel
      &\global\@eqcnt\@ne \hfil$\displaystyle{{}##{}}$\hfil
      &\global\@eqcnt\tw@ $\displaystyle{##}$\hfil\tabskip\@centering
      &\global\@eqcnt\thr@@ \hb@xt@\z@\bgroup\hss##\egroup
        \tabskip\z@skip
      \cr}
\makeatother
\makeatletter
\def\varin{\mathrel{\mathpalette\@varin\relax}}
\def\@varin#1{%
   \hbox{\setbox\z@\hbox{\m@th$#1\cup$}%
       \def\reserved@a{bold}%
       \dimen@\ifx\reserved@a\math@version .3\else .2\fi\p@
       \kern.5\wd\z@\kern-\dimen@
       \vrule\@width2\dimen@\@height1.08\ht\z@\@depth\z@
       \kern-\dimen@\kern-.5\wd\z@
       \box\z@}}
\makeatother

\begin{document}
 


\title [Atiyah-Patodi-Singer condition] 
{Atiyah-Patodi-Singer boundary condition 
and a splitting formula of a spectral flow}



\author{Kenro Furutani}
\address{Kenro Furutani \endgraf 
Department of Mathematics \endgraf 
Faculty of Science and Technology \endgraf 
Science University of Tokyo \endgraf 
2641 Noda, Chiba (278-8510)\endgraf 
Japan \endgraf}
\email{furutani@ma.noda.sut.ac.jp}

\bigskip
\bigskip

\begin{abstract}
We describe a relation between Atiyah-Patodi-Singer boundary 
condition and a
global elliptic boundary condition which naturally
appears in formulating a splitting formula for a spectral flow,
when we decompose the manifold into two 
components. Then we give a variant of the splitting formula
with the H\"ormander index as a correction term.
\end{abstract}

\maketitle
\tableofcontents
\thispagestyle{empty}

\section{Introduction}
In the paper \cite{FO} we formulated and 
proved a splitting formula of a
spectral flow for a continuous family of first order 
selfadjoint elliptic differential operators 
$\{A_t\}_{t\in [0,1]}$ defined 
on a closed manifold $M$. This is an addition formula of a spectral flow 
when we decompose a manifold 
into two components along a hypersurface 
$\Sigma$, $M=M_-\cup_{\Sigma} M_+$, $\partial M_{\pm}=\Sigma$.
For such a family $\{A_t\}$ considered on the whole closed manifold $M$, 
an integer, called spectral flow, is well defined,
and to ``observe'' this quantity we cut the manifold along a
hypersurface, then we can ``observe'' from the hypersurface
a quantity ``Maslov index''
which is a curve of boundary data of solutions of operators. This
quantity can be understood as the spectral flow. This is just a spectral
flow formula(\cite{Yo}, \cite{Ni}) where manifolds need not be
separated
into two components. If the manifold is separated into two parts
by the hypersurface,
then we will have two Maslov indexes which together give the
whole information of the spectral flow.
For this observation we must make clear which family we are
observing, i.e.,
to get a family of
selfadjoint Fredholm 
extensions we must impose 
a suitable elliptic boundary condition 
on the family $\{A_t\}$ 
when we restrict the operators on each component $M_{\pm}$.
This condition appears in a natural way 
in our formulation to write
down the splitting formula and reflects the influence
from one side to other side.
On the other hand Atiyah-Patodi-Singer boundary condition
is described based on the boundary data. Nevertheless
these two are relating each other in the case of
the operators of the product form near the hypersurface.

The purpose of this paper is to describe a relation
between our global elliptic boundary condition
and Atiyah-Patodi-Singer boundary condition.
There are several such splitting formulas (\cite{CLM2}, \cite{DK},
\cite{KL},\cite{Ta},\cite{Ni}) 
and the boundary condition treated there is mostly
Atiyah-Patodi-Singer boundary condition. So the result in this paper
will give us a sight of our global elliptic
boundary condition for operators satisfying suitable
analytic assumptions.

We follow the theory of symplectic Hilbert spaces, especially 
of the Fredholm-Lagrangian-Grassmannian 
and the Maslov index in the infinite dimension which were
discussed in the papers \cite{Fu} and \cite{FO} precisely(also see
\cite{Go}, \cite{CLM1} and \cite{RS}).

In $\S 2$ we explain a global elliptic boundary condition appearing in
the splitting formula for a spectral flow and state a relation between
it and 
Atiyah-Patodi-Singer boundary condition in terms of the
Fredholm-Lagrangian-Grassmannian.

In $\S 3$, first we recall a $L_2$-reduction theorem \cite{Fu} and by applying
this we give a proof of Theorem (\ref{thm-1}).

In $\S 4$ as an application of Theorem (\ref{thm-1}) we give
a variant of a splitting formula of a spectral flow with H\"ormander
index as a correction term.



\section{A global elliptic boundary condition}

In this section we explain an elliptic boundary condition we
introduced in \cite{FO}.

Let $A$ be a first order selfadjoint elliptic differential 
operator defined on a
real vector bundle $\mathbb{E}$ on a closed manifold $M$.
Let $\Sigma$ be a hypersurface of $M$ along which $M$ is separated
into two components $M_{\pm}$, $M=M_-\cup_{\Sigma}M_+$, $\partial
M_{\pm}=\Sigma$,
and we denote the first order Sobolev space on $M$ (resp. $M_{\pm}$)
taking values in the real vector bundle $\mathbb{E}$ by $H^1(M,\mathbb{E})$ 
(resp. $H^1(M_{\pm},\mathbb{E}_{|M_{\pm}}$)).  For the subspace in 
$H^1(M_{\pm},\mathbb{E}_{|M_{\pm}})$ with vanishing boundary values
we denote it by $H^1_0(M_{\pm},\mathbb{E}_{|M_{\pm}})$.
These are the domains of the minimal closed extensions of the
operators $A$ considered on
$C^{\infty}_0(M_{\pm}\backslash\Sigma,\mathbb{E}_{|M_{\pm}\backslash\Sigma})$ and
we denote them by $\mathfrak{D}_{\bf m}^{\pm}$ =
$H^1_0(M_{\pm},\mathbb{E}_{|M_{\pm}})$. 
Then we denote by $\mathcal{A}_{\pm}^{\,\,\,*}$
their adjoint operators considered on $M_{\pm}$ and by
$\mathfrak{D}_{\bf M}^{\pm}$
their domains of definitions, i.e., $f\in\mathfrak{D}^{\pm}_{\bf M}$
if $f\in L_2(M_{\pm},\mathbb{E}_{|M_{\pm}})$ 
and $A(f)\in L_2(M_{\pm},\mathbb{E}_{|M_{\pm}})$ in the distribution sense. 

We must put two assumptions {\bf (a1)} and
{\bf (a2)} on the selfadjoint elliptic operator $A$:\\

\qquad {\bf (a1)}: $A$ satisfies the {\it unique continuation property} with respect
   to the hypersurface $\Sigma$, that is,
\begin{equation}
\Ker(\mathcal{A}_{\pm}^{\,\,\,*})\cap 
\mathfrak{D}^{\pm}_{\bf m}=\{0\}.
\end{equation}

\qquad {\bf (a2)}: On a tubular 
neighborhood $\mathcal{N}\cong (-1,1)\times \Sigma$
the operator $A$ is of the product form, that is,
\begin{equation}
A= \sigma\Bigr(\frac{\partial}{\partial u}+ B_0\Bigr),
\end{equation}
where $\sigma$ is a bundle map on $\mathbb{E}_{|\mathcal{N}}$ which
does not depend on the normal variable $u\in (-1,1)$, the operator
$B_0$ is a selfadjoint elliptic operator on $\Sigma$ and also
does not depend on the normal variable $u$.

We identify $\mathcal{N}\cap M_- \cong(-1,0]\times \Sigma$ and
$\mathcal{N}\cap M_+ \cong[0,1)\times \Sigma$.

The selfadjointness of the operator $A$ implies that
\begin{align}
&\sigma^2=-\Id,~^t\sigma=-\sigma\label{property-1}\\
&\sigma \circ B_0 +  B_0\circ \sigma =0,\label{property-2}
\end{align}
where the transpose $^t\sigma$ is taken with respect 
to a suitable metric on the vector bundle $\mathbb{E}$. 
We denote the inner product on the $L_2$ sections of $\mathbb{E}$
by $<\cdot,\cdot>$. Then $\sigma$ defines an almost complex structure
and a compatible symplectic structure
on $L_2(\Sigma, \mathbb{E}_{|\Sigma)})=L_2(\Sigma)$.

Let $\mathfrak{D}_0$ be a subspace in $H^1(M_{-},\mathbb{E}_{|M_{-}})$
defined by
\[
\mathfrak{D}_0
=\bigr\{ f\in H^1(M_{-},\mathbb{E}_{|M_{-}})\,\bigr|\, ^{\exists}
\tilde{f}\in H^1(M,\mathbb{E}) ~\text{such that}~\tilde{f}_{\,\,\,|M_-}=f,~\text{and}~A(f)=0 ~\text{on} ~M_+\bigr\}.
\]

Hereafter we will state properties only for the case $M_-$,
but shall use the corresponding results for $M_+$, if necessary.

We denote the restriction of $\mathcal{A}_{-}^{\,\,\,*}$ to $\mathfrak{D}_0$ 
by $A_{\mathfrak{D}_0}$,
then under the two assumptions {\bf (a1)}
and {\bf (a2)} we have
\begin{prop}
The operator $A_{\mathfrak{D}_0}$ 
satisfies the inequality:
\begin{equation}
\Vert f \Vert_1\leq c(\Vert A_{\mathfrak{D}_0}(f)\Vert +\Vert
f\Vert ), ~\forall f\in \mathfrak{D}_0
\end{equation}
with a positive constant $c>0$,
where $\Vert f\Vert_1$ denotes the first order Sobolev norm.
And so $A_{\mathfrak{D}_0}$ is selfadjoint and has compact resolvents. 
\end{prop}

This property is basic to state Theorem (\ref{thm-1}) below.
We already made use of this property 
in our paper (\cite{FO}). The proof is given upon 
$L_2$-reduction theorem of the Maslov index in the infinite
dimension. In the next section
we explain some part of a proof of this proposition, and together
with the help of
Rellich's Theorem we prove Theorem (\ref{thm-1}).

\begin{rem}
If $A$ is invertible on $H^1(M,\mathbb{E})$, then
the same holds for $A_{\mathfrak{D}_0}$ as in the above proposition
without assumptions {\bf (a1)} and {\bf (a2)}, because
we have 
\[
\Vert f\Vert_{H^1(M_-,\mathbb{E}_{M_-})}\le C(\Vert
f\Vert_{H^1(M,\mathbb{E})}+\Vert f\Vert_{L_2(M)})=
C(\Vert
f\Vert_{H^1(M_-,\mathbb{E}_{|M_-})}+\Vert f\Vert_{L_2(M)})
\le C'\Vert
f\Vert_{H^1(M_-,\mathbb{E}_{|M_-})}
\]
for $f\in H^1(M,\mathbb{E})$ satisfying $A(f)=0 $ on $M_+$.
\end{rem}

{\it Let $\{\ell_k\}_{k\in\mathbb{Z}\backslash\{0\}}$ 
be the eigenvalues of the operator $B_0$
and we denote corresponding orthonormal 
eigensections by $\{\varphi_k\}$.}
{}From the properties (\ref{property-1}) and (\ref{property-2}) 
we have  
$\ell_k=-\ell_{-k}> 0$ for $k=N_0+1,N_0+2,\cdots $ with 
$N_0=1/2\cdot\dim\Ker(B_0)$ (see Remark (\ref{zero-mode}) below)
and $\ell_k=0$ for $0<|k|\le N_0$.

For a section $\varphi$ on $\Sigma$ let 
\[
\varphi=\sum\limits_{k
  \in\mathbb{Z}\backslash\{0\}}
a_k\varphi_k
\]
be the eigensection-expansion, then 
the Sobolev space 
$H^s(\Sigma,\mathbb{E}_{|\Sigma})$ of order $s\in\mathbb{R}$ 
on $\Sigma$ is characterized as
\[
H^s(\Sigma,\mathbb{E}_{|\Sigma})=\bigr\{\,
\varphi=\sum\limits_{k
  \in\mathbb{Z}\backslash\{0\}}
a_k\varphi_k\,\bigr|\,
\sum\limits_{k\in\mathbb{Z}\backslash\{0\}}
|a_k|^2|\ell_k|^{2s}<\infty\bigr\}.
\]

Let $\mathfrak{D}_{APS}^{\,\,0}$ 
be a subspace in $H^1(M_-,\mathbb{E}_{|M_-})$
such that
\begin{align}\label{aps-0}
&\mathfrak{D}_{APS}^{\,\,0}\\
& = \bigr\{ 
f\in H^1(M_-,\mathbb{E}_{|M_-})~\bigr|~\text{if}~f_{\,|\Sigma}=\sum\limits_{k
  \in\mathbb{Z}\backslash\{0\}}
a_k\varphi_k,~\text{then}~a_{-k}=0 ~\text{for}~k=1,2,3,\cdots
 \bigr\},\notag
\end{align}
and we denote the restriction of $\mathcal{A}_-^{\,\,\,*}$ to
$\mathfrak{D}_{APS}^{\,\,0}$
by $A_{\mathfrak{D}_{APS}^{\,\,0}}$. This non-local 
boundary condition (\ref{aps-0}) is called
``Atiyah-Patodi-Singer boundary condition''. Then
\begin{prop}{\bf(\cite{APS})}
The operator $A_{\mathfrak{D}_{APS}^{\,\,0}}$
defined on $\mathfrak{D}_{APS}^{\,\,0}$ 
is a selfadjoint operator 
with compact resolvents.
\end{prop}

Now let $\beta^-=\mathfrak{D}^-_{\bf M}/\mathfrak{D}^-_{\bf m}$ be the
space of boundary values.  Here the maximal domain 
$\mathfrak{D}^-_{\bf M}$ is
equipped
with the norm $\Vert\cdot\Vert^{G}$ given by the graph inner product:
\[
<f,g>^{G}=<f,g>+<\mathcal{A}_-^{\,\,\,*}(f),\mathcal{A}_-^{\,\,\,*}(g)>.
\]
The space $\beta^-$ has a structure of a symplectic Hilbert space with
the symplectic form
\begin{equation}\label{symplectic-form}
\omega([f],[g])=<\mathcal{A}_-^{\,\,\,*}(f),g> -
<f,\mathcal{A}_-^{\,\,\,*}(g)>,~f,~g\in \mathfrak{D}^-_M,
\end{equation}
and is realized in the distribution space on $\Sigma$:
\begin{align*}
&\beta^- =\bigr\{f\in H^{-1/2}(\Sigma,\mathbb{E}_{|\Sigma})\,\bigr|\,
f=\sum\limits_{k\in\mathbb{Z}\backslash\{0\}}c_k\varphi_k,~\text{with}~
\sum\limits_{k>0}|c_k|^2\ell_k<\infty ~\text{and}\\
&\qquad\qquad\qquad\qquad\qquad\qquad\qquad\qquad
\sum\limits_{k<0}|c_k|^2|\ell_k|^{-1}<\infty
\bigr\}\\
&=\bigr\{f\in H^{1/2}(\Sigma,\mathbb{E}_{|\Sigma})\,\bigr|\,
f=\sum\limits_{k>0}c_k\varphi_k,~\text{with}~
\sum\limits_{k>0}|c_k|^2\ell_k<\infty\bigr\}\\
&+\bigr\{f\in H^{-1/2}(\Sigma,\mathbb{E}_{|\Sigma})\,\bigr|\,
f=\sum\limits_{k<0}c_k\varphi_k,~\text{with}~
\sum\limits_{k<0}|c_k|^2|\ell_k|^{-1}<\infty\bigr\}\\
&=\theta^-_{+}+\theta^-_{-}.
\end{align*}

For the determination of the space $\beta^{-}$ 
see \cite{Ho1} and \cite{APS}. 

\begin{rem}\label{zero-mode}
By the conditions {\bf (a1)} and {\bf (a2)} 
we know 
that $\Ker(B_0)$
is a finite dimensional symplectic subspace of $\beta^-$, so
that we choose eingensections $\{\varphi_k\}$ for $0<|k|\leq N_0$ in 
such a way that the subspaces spanned by 
$\{\varphi_k\}_{-N_0\leq k <0}$
and $\{\varphi_k\}_{0<k\leq N_0}$ are mutually transversal
Lagrangian subspaces in $\Ker(B_0)$.
\end{rem}

\begin{rem}\label{synplectic form by sigma}
Of course for smooth sections
of $\mathbb{E}_{|\Sigma}$(and also
for $L_2$-sections) the symplectic structure
defined by $\sigma$ coincides with $\omega$ 
defined in (\ref{symplectic-form}).
\end{rem}

Let $\gamma^- : \mathfrak{D}^-_{\bf M}
\rightarrow \beta^-$ be the projection map, then the image 
$\gamma^-(\Ker (\mathcal{A}_-^{\,\,\,*}))$ 
is a Lagrangian subspace and the pairs 
$\bigr(\gamma^-(\Ker (\mathcal{A}_-^{\,\,\,*})),
\gamma^-(\mathfrak{D}_{APS}^{\,\,0})\bigr)$
and 
$\bigr(
\gamma^-(\Ker (\mathcal{A}_-^{\,\,\,*})),\gamma^-(\mathfrak{D}_0)
\bigr)$ 
are Fredholm pairs. The spaces 
$\gamma^{-}\bigr(\Ker(\mathcal{A}_{-}^{\,\,\,*})\bigr)$
is called {\it Cauchy data space}. 
Further, for each Lagrangian subspace $\lambda\subset\beta^-$
the operator ${\mathcal{A}_{-}^{\,\,\,*\,\,}}_{\bigr|(\gamma^-)^{-1}(\lambda)}$
is a selfadjoint realization and if 
$\bigr(\lambda, \gamma^{-}(\Ker(\mathcal{A}_-^{\,\,\,*}))\bigr)$
is a Fredholm pair, then 
${\mathcal{A}_-^{\,\,\,*}}_{\bigr|(\gamma^-)^{-1}(\lambda)}$
is a selfadjoint Fredholm operator.

\begin{rem}
Note that the Lagrangian property of the Cauchy data space 
$\gamma^{-}\bigr(\Ker(\mathcal{A}_{-}^{\,\,\,*})\bigr)$
would not be trivial. To prove 
this property we rely on the existence
of at least one selfadjoint realization
of $A_{\,|\mathfrak{D}_{\bf m}^{-}}$ 
(= restriction of $\mathcal{A}_{-}^{\,\,\,*}$ to a suitable subspace
in $\mathfrak{D}_{\bf M}^{-}$) with compact resolvents
(or a unbounded selfadjoint Fredholm extension). 
For our case this realization is 
given by the operator $A_{\mathfrak{D}_{APS}^{\,\,0}}$, and 
for which proof we use the assumption {\bf (a2)}. 
\end{rem}

We denote by $\mathcal{F}\Lambda_{\lambda}(\beta^-)$
the space of Lagrangian subspaces $\mu$ in $\beta^-$
such that 
$(\mu,\lambda)$ is a Fredholm pair:
\[
\mathcal{F}\Lambda_{\lambda}(\beta^-)
=\bigr\{\mu\subset\beta^-\,\bigr|\, \mu~\text{is a Lagrangian subspace
  and}~(\mu,\lambda)~
\text{is a Fredholm pair}
\bigr\}.
\]
 
Now we have 
\begin{thm}\label{thm-1}
\begin{equation}
\mathcal{F}\Lambda_{\gamma^-(\mathfrak{D}_0)}(\beta^-)
=\mathcal{F}\Lambda_{\gamma^-(\mathfrak{D}_{APS}^{\,\,0})}(\beta^-),
\end{equation}
more precisely the orthogonal projection operators onto the subspace
$\gamma^-(\mathfrak{D}_0)$ and that onto the subspace 
$\gamma^-(\mathfrak{D}_{APS}^{\,\,0})=\theta^-_+$ 
differ by a compact operator.
\end{thm}

We prove this in the next section.

\begin{rem}\label{flg}
Let $H$ be a symplectic Hilbert space
and we regard $H$ as a complexification of a Lagrangian subspace
$\lambda$. Then a Lagrangian subspace $\mu$ is of the form  
$U(\lambda)=\mu$  with a unitary operator $U$ of the form 
$\Id$+ {\it compact operator}, 
then
$\mathcal{F}\Lambda_{\lambda}(H)=\mathcal{F}\Lambda_{\mu}(H)$.
Also this property is equivalent to the condition that
the difference of the orthogonal projection operators
onto the Lagrangian subspaces $\lambda$ and $\mu$
is compact. 
For such two Lagrangian subspaces $\lambda$ and $\mu$ and
two arbitrary Lagrangian subspaces $\nu_0$ and $\nu_1$ in
$\mathcal{F}\Lambda_{\lambda}(H)=\mathcal{F}\Lambda_{\mu}(H)$
we have a well-defined integer $\sigma(\nu_0,\nu_1;\lambda,\mu)$,
called H\"ormander index. This is the difference
of the Maslov indexes 
\[
\sigma_H\bigr(\nu_0,\nu_1;\lambda,\mu\bigr)=
\Mas\bigr(\{c(t)\},\lambda\bigr)-\Mas\bigr(\{c(t)\},\mu\bigr),
\]
where the path $\{c(t)\}$ is in $\mathcal{F}\Lambda_{\lambda}(H)$
connecting $\nu_0$ and $\nu_1$ and the difference does not
depend on any such paths.
Here the Maslov index $\Mas\bigr(\{c_t\},\lambda\bigr)$
is , in a sense, the intersection number with the ``Maslov cycle''
$\mathfrak{M}_{\lambda}
=\bigr\{\mu\in\mathcal{F}\Lambda_{\lambda}(H)\,\bigr|\,
\mu\cap\lambda\not=\{0\}\bigr\}$(\cite{Ho2}, \cite{Fu}).
\end{rem}

\section{Symplectic reduction theorem}

In this section after recalling 
a symplectic reduction theorem (\cite{Fu})
we prove Theorem (\ref{thm-1}).

Let $(\mathcal{B}, \omega_{\mathcal{B}})$ and $(L,\omega_L)$ 
be two symplectic Hilbert spaces ($\omega_{\mathcal{B}}$ is the
symplectic form and so on)
with decompositions by Lagrangian subspaces 
$\theta_-$, $\theta_+$, $L_-$ and $L_+$
({\it Polarized symplectic Hilbert space}):
\begin{equation}\label{decompo-1}
\mathcal{B}=\theta_- +\theta_+,~L=L_- + L_+.
\end{equation}
We assume that there are continuous injective maps 
${\bf i}_{+}:L_+ \to \theta_+$ and ${\bf i}_-:\theta_- \to L_-$ 
having dense images such that
\begin{equation}\label{decompo-2}
\omega_{\mathcal{B}}({\bf i}_+(a), x)=
\omega_L(a,{\bf i}_-(x))~\text{for any} ~x\in \theta_-~\text{and}~a\in L_+.
\end{equation}
Then 
\begin{prop}{\bf (\cite{Fu})}
There is a continuous map 
$\tau :\mathcal{F}\Lambda_{\theta_-}(\mathcal{B}) \to
\mathcal{F}\Lambda_{L_-}(L)$ such that
for any continuous curve $\{c(t)\}_{t\in[0,1]}$ in
$\mathcal{F}\Lambda_{\theta_-}(\mathcal{B})$
\[
\Mas\bigr(\{c(t)\},\theta_-\bigr)= \Mas\bigr(\{\tau(c(t))\},L_-\bigr).
\]
\end{prop}

The map $\tau$ is defined in the following way:
\[
\tau(\nu)=\bigr\{ b+a \in L=L_+ + L_-\,\bigr|\,
{^{\exists}{x}}\in \theta_- ~\text{such that}~
{\bf i}_+(b)+ x\in\nu ~\text{and}~a={\bf i}_-(x)\bigr\}.
\]

For any decomposition of $\theta_-=F+F'$ by closed subspaces
$F$($\dim F<+\infty$) and $F'$,
we can decompose $L_+$ by closed subspaces in such a way that
$L_+=G+G'$ with $\dim G =\dim F$ and $G+{\bf i}_-(F)$ is a symplectic
subspace in $L$. Also in this case the subspace $F+{\bf i}_+(G)$
is a symplectic subspace in $\mathcal{B}$.
Moreover the subspaces $F+\overline{{\bf i}_+(G')}$ and
$G+\overline{{\bf i}_-(F')}$ are Lagrangian subspaces.

Then by replacing $\theta_-$ with $F'+{\bf i}_+(G)$,
$\theta_+$ with $F+\overline{{\bf i}_+(G')}$,
$L_-$ with $G+\overline{{\bf i}_+(F')}$
and $L_+$ with ${\bf i}_-(F)+ G'$ and also by replacing the maps
${\bf i}_{\pm}$ in an obvious way 
we have a similar situation as in
(\ref{decompo-1}) and (\ref{decompo-2}). 
We shall denote these new maps by
$\tilde{{\bf i}}_{\pm}$, although
the resulting maps $\tau$ between Fredholm-Lagrangian-Grassmannians
$\mathcal{F}\Lambda_{\theta_-}(\mathcal{B})=
\mathcal{F}\Lambda_{F'+{\bf i}_+(G)}(\mathcal{B})$
and
$\mathcal{F}\Lambda_{L_-}(L)=
\mathcal{F}\Lambda_{G+\overline{{\bf i}_-(F')}}(L)$
coincides. 

Note that the arguments above are guaranteed that the
spaces $\mathcal{B}$ and $L$ are Hilbert spaces 
(see \cite{KS} for symplectic Banach spaces).

We apply this proposition to the case
$\mathcal{B}=\beta^+=\theta^+_{+} + \theta^+_-$ and
$L=\beta^- = \theta^-_{+} +\theta^-_-$.
Note that the space $\beta^+$ is defined as follows:
\begin{align*}
&\beta^+ =\bigr\{f\in H^{-1/2}(\Sigma,\mathbb{E}_{|\Sigma})\,|\,
f=\sum\limits_{k\in\mathbb{Z}\backslash\{0\}}c_k\varphi_k,~\text{with}~
\sum\limits_{k>0}|c_k|^2\ell_k^{-1}<\infty ~\text{and}\\
&\qquad\qquad\qquad\qquad\qquad\qquad\qquad\qquad
\sum\limits_{k<0}|c_k|^2|\ell_k|<\infty
\bigr\}\\
&=\bigr\{f\in H^{-1/2}(\Sigma,\mathbb{E}_{|\Sigma})\,\bigr|\,
f=\sum\limits_{k>0}c_k\varphi_k,~\text{with}~
\sum\limits_{k>0}|c_k|^2\ell_k^{-1}<\infty\bigr\}\\
&+\bigr\{f\in H^{1/2}(\Sigma,\mathbb{E}_{|\Sigma})\,\bigr|\,
f=\sum\limits_{k<0}c_k\varphi_k,~\text{with}~
\sum\limits_{k<0}|c_k|^2|\ell_k|<\infty\bigr\}\\
&=\theta^+_{+}+\theta^+_{-}.
\end{align*}
 
The maps ${\bf i}_{\pm}$ here are given by inclusion maps. 

Since ($\gamma^+(\Ker(\mathcal{A}_+^{\,\,\,*})), \theta^+_-$)
is a Fredholm pair, we can find a finite dimensional subspace $F$
in $\theta^+_{-}$ 
and a corresponding finite dimensional subspace
$G$ in $\theta^-_{+}$ such that
we have  decompositions
\[
\theta^+_{-}=F+F',\quad \theta_{+}^-=G+G'
\]
with suitable closed subspaces $F'$ and $G'$ and that
\begin{center}
$F'+ {\bf i}_{+}(G)$ and $\gamma^+(\Ker(\mathcal{A}_+^{\,\,\,*}))$ 
are transversal.
\end{center}

When we put $F'+ {\bf i}_{+}(G)=\lambda_-$ and
$F + \overline{{\bf i}_+(G')}=\lambda_+$
we have a decomposition $\beta^+=\lambda_+ +\lambda_-$ 
with Lagrangian subspaces $\lambda_{\pm}$ and the Cauchy data space
is expressed as a graph of a continuous 
map $\mathfrak{K}: \lambda_+ \to
\lambda_-$. Then for such a Lagrangian subspace we have
$\tau(\gamma^+(\Ker(\mathcal{A}_+^{\,\,\,*})))$ 
is the graph of the map $\tilde{{\bf i}}_-\circ \mathfrak{K}\circ \tilde{{\bf i}}_+$.
Note here the maps $\tilde{{\bf i}}_{\pm}$ should be defined in a
suitable way according to the choices of the subspace $F$ and $G$
(for example, ${\tilde{\bf i}_+}$ is defined as
${\bf i}_+$ on $F'$ and 
${\bf i}_-^{\,\,-1}$ on ${\bf i}_-(G)$).

Now the original maps ${\bf i}_{+}:\theta^-_+\to \theta^+_+$  and
${\bf i}_{-}:\theta^+_-\to \theta^-_-$ 
are compact operators by Rellich's Theorem
and so the new maps $\tilde{{\bf i}}_{\pm}$ are also compact. 

Let us denote the orthogonal projection operator to a closed subspace
$E$ by $\mathcal{P}_E$.
Then the difference 
\begin{equation}\label{difference-1}
\mathcal{P}_{{\bf i}_+(F)+G'}
-\mathcal{P}_{\tau(\Ker(\mathcal{A}_+^{\,\,\,*}))}
\end{equation}
is a compact operator and the difference
\begin{equation}\label{difference-2}
\mathcal{P}_{{\bf i}_+(F)+G'}-\mathcal{P}_{\theta^-_+}
\end{equation}
is a finite rank operator.

By the definition of the map $\tau$ 
we have $\gamma^-(\mathfrak{D}_0)=
\tau(\gamma^+(\Ker(\mathcal{A}_+^{\,\,\,*})))$, and
(\ref{difference-1}) and (\ref{difference-2}) imply
that the difference of the orthogonal
projection operators onto the subspaces
$\gamma^-(\mathfrak{D}_{APS}^{\,\,0})$ = $\theta^-_+$ 
and $\gamma^-(\mathfrak{D}_0)$ is a compact operator.

{\it So this gives us a proof of Theorem (\ref{thm-1}).}

\section{Cauchy data spaces and H\"ormander index}

Let $L_2(\Sigma)= L_+ +L_-$ be the polarization 
by $L_{\pm}$, where $L_+$ is the $L_2$-completion
of the space spanned by $\{\varphi_k\}_{k>0}$ 
and $L_-$ is the $L_2$-completion
of the space spanned by $\{\varphi_k\}_{k<0}$.
Then by applying above arguments to the two pairs 
$\bigr(\beta^+=\beta^+_- +\beta^+_+, L_2(\Sigma)=L_- +L_+\bigr)$
and $\bigr(\beta^-=\beta^-_- +\beta^-_+, L_2(\Sigma)=L_- +L_+\bigr)$
of polarized symplectic Hilbert spaces we have
four Lagrangian subspaces
\[
\gamma^{\pm}\bigr(\Ker(\mathcal{A}_{\pm}^{\,\,\,*})\bigr)
\cap L_2(\Sigma), 
~L_{\pm}
\]
of $L_2(\Sigma)$ which satisfy following properties {\bf (h1)} and
{\bf (h2)}:
\begin{align*}
&{\bf (h1)}:\gamma^{\pm}\bigr(\Ker(\mathcal{A}_{\pm}^{\,\,\,*})\bigr)
\cap L_2(\Sigma)~\text{and}~
L_{\mp}~\text{are Fredholm pairs},\\
&{\bf (h2)}:\gamma^{\pm}\bigr(\Ker(\mathcal{A}_{\pm}^{\,\,\,*})\bigr)
\cap L_2(\Sigma)
=U_{\pm}(L_{\pm}),~\text{where}~U_{\pm}
~\text{are unitary operators of}\\
&\text{the form Id +{\it compact operator}}.
\end{align*}
Here we identify
$L_2(\Sigma)\cong L_+\otimes\mathbb{C}$.

Now we can define the H\"ormander index
\begin{equation}\label{asymmetry}
\sigma_H\bigr(\gamma^+(\Ker(\mathcal{A}_{+}^{\,\,\,*}))
\cap L_2(\Sigma), L_{+};
\gamma^-(\Ker(\mathcal{A}_{-}^{\,\,\,*}))\cap L_2(\Sigma), 
L_{-}\bigr)
\end{equation}
of these four Lagrangian subspaces. Then 
its absolute value will express an {\it asymmetry} of solution spaces
of the operator $A$ under the decomposition of $M$
along a hypersurface $\Sigma$. So, if there is a {\it symmetry}
among these four Lagrangian subspaces, the value must vanish. In fact

\begin{prop}
Assume that 
$\sigma\bigr(\gamma^+(\Ker(\mathcal{A}_+^{\,\,\,*}))\cap L_2(\Sigma)\bigr)$
$=\gamma^-\bigr(\Ker(\mathcal{A}_-^{\,\,\,*})\bigr)\cap L_2(\Sigma)$, then
the H\"ormander index of these four Lagrangian subspaces vanishes:
\[
\sigma_H\bigr(\gamma^+(\Ker(\mathcal{A}_{+}^{\,\,\,*}))
\cap L_2(\Sigma), L_{+};
\gamma^-(\Ker(\mathcal{A}_{-}^{\,\,\,*}))\cap L_2(\Sigma), 
L_{-}\bigr)=0.
\]
\end{prop}

\begin{proof}
First we assume that 
$\gamma^+(\Ker(\mathcal{A}_+^{\,\,\,*}))\cap L_2(\Sigma)$ and $L_-$
are transversal. Then
the space $\gamma^+(\Ker(\mathcal{A}_+^{\,\,\,*}))\cap L_2(\Sigma)$ 
is written as a graph of a compact operator 
$T:L_+ \to L_-$ such that $\sigma\circ T$
is a selfadjoint operator on $L_+$ and
the space 
$\gamma^-(\Ker(\mathcal{A}_-^{\,\,\,*}))\cap L_2(\Sigma)$
is also written as a graph of the map 
$-\sigma\circ T\circ \sigma$. These imply that
the curve of Lagrangian subspaces given by the graphs of 
$\bigr\{-t \cdot \sigma \circ T\circ \sigma
\bigr\}_{0\leq t \leq 1}$ 
is always transversal
to both of 
$\gamma^+(\Ker(\mathcal{A}_+^{\,\,\,*}))\cap L_2(\Sigma)$
and $L_+$. This curve is connecting 
$\gamma^-(\Ker(\mathcal{A}_-^{\,\,\,*}))\cap L_2(\Sigma)$
and $L_-$.
Hence we have
\[
\sigma_H\bigr(\gamma^-(\Ker(\mathcal{A}_{-}^{\,\,\,*}))
\cap L_2(\Sigma), L_{-};
\gamma^+(\Ker(\mathcal{A}_{+}^{\,\,\,*}))\cap L_2(\Sigma), 
L_{+}\bigr)=0.
\]

If $\gamma^+(\Ker(\mathcal{A}_+^{\,\,\,*}))\cap L_2(\Sigma)$
and $L_-$ are not transversal, then we decompose
the Lagrangian subspace
$\gamma^+(\Ker(\mathcal{A}_+^{\,\,\,*}))\cap L_2(\Sigma)$
into the orthogonal sum
\[
\gamma^+(\Ker(\mathcal{A}_+^{\,\,\,*}))\cap L_2(\Sigma)
=\ell_0+\nu,
\]
where
$
\ell_0=\bigr(\gamma^+(\Ker(\mathcal{A}_+^{\,\,\,*}))
\cap L_2(\Sigma)\bigr)\cap L_- 
$
and
$\nu$ is the orthogonal complement of $\ell_0$ in 
$\gamma^+(\Ker(\mathcal{A}_+^{\,\,\,*}))\cap L_2(\Sigma)$.
Also we decompose $L_-=\ell_0 + (L_-\cap \ell_0^{\,\,\perp})=
\ell_0+\ell_-$
and $L_+=\sigma(\ell_0)+(L_+\cap\sigma(\ell_0)^{\perp})=
\sigma(\ell_0)+\ell_+$.
Now we have
\begin{align*}
&\sigma_H\bigr(\gamma^-(\Ker(\mathcal{A}_{-}^{\,\,\,*}))
\cap L_2(\Sigma), L_{-};
\gamma^+(\Ker(\mathcal{A}_{+}^{\,\,\,*}))\cap L_2(\Sigma), 
L_{+}\bigr)\\
&=\sigma_H\bigr(\sigma(\ell_0),\ell_0 ;\ell_0,\sigma(\ell_0)\bigr)
+\sigma_H\bigr(\sigma(\nu),\ell_- ;\nu,\ell_+\bigr)=0,
\end{align*}
by applying the first arguments to the second term.

Note that $\gamma^-(\Ker(\mathcal{A}_{-}^{\,\,\,*}))
\cap L_2(\Sigma)=\sigma(\ell_0)+\sigma(\nu)$ 
is an orthogonal decomposition and the vanishing of the first term
follows from a skew-symmetric property of the H\"ormander
index.
\end{proof}

\section{A splitting formula of a spectral flow}

First we state a splitting formula for a spectral flow when
we decompose a manifold into two components.
Then we give another form of it by replacing the boundary
condition with Atiyah-Patodi-Singer condition.

Let $\{C_t\}_{t\in[0,1]}$ be a continuous family of 
symmetric bundle maps
of $\mathbb{E}$ and we assume that each of the operator in 
the family $\{A+C_t\}$ satisfies the
conditions {\bf (a1')} and {\bf (a2)} where {\bf (a1')} is:

\qquad {\bf (a1')}: There exists an $\epsilon_0>0$ 
such that for any $|s|< \epsilon_0$ and any $t\in [0,1]$
the operators $A+C_t+s$ satisfy 
the {\it unique continuation property}
with respect to the hypersurface $\Sigma$ : 
\begin{equation}\label{strong-unique-cont}
\Ker(\mathcal{A}_{\pm}^{\,\,\,*}+C_t+s)\cap \mathfrak{D}_{\bf m}^{\pm}=\{0\}.
\end{equation}
 
Here $C_t$ is regarded as a bounded selfadjoint 
operator on $L_2(M,\mathbb{E})$. 

Now we have continuous families of Cauchy data spaces
$\gamma^{\pm}(\Ker(\mathcal{A}_{\pm}^{\,\,\,*}+C_t))$ 
($C_t$ should be considered 
as acting on the space $\mathfrak{D}_{\bf M}^{\pm}$ respectively, 
and both of which are invariant under this action). The splitting
formula is stated as follows:
\begin{thm}{\bf(\cite{FO})}\label{thm-2}
\begin{equation}
\Sf\bigr(\{A+C_t\}\bigr)=\Sf\bigr(\{A_{\mathfrak{D}_0}+C_t\}\bigr)+
\Sf\bigr(\{A_{\mathfrak{D}_1}+C_t\}\bigr),
\end{equation}
where
\[
\mathfrak{D}_0
=\bigr
\{f\in 
H^1(M_-,\mathbb{E}_{|M_-})\,
\bigr|\,{^\exists \tilde{f}}\in H^1(M,\mathbb{E)}~\text{such that}
~\tilde{f}_{\,\,|M_-}=f 
~\text{and}~(A+C_0)(\tilde{f})=0 ~\text{on}~M_+\bigr\}
\]
and
\[
\mathfrak{D}_1
=\bigr
\{g\in 
H^1(M_+,\mathbb{E}_{|M_+})\,
\bigr|\,{^\exists \tilde{g}}\in H^1(M,\mathbb{E)}~\text{such that}
~\tilde{g}_{|\,M_+}=g 
~\text{and}~(A+C_1)(\tilde{g})=0 ~\text{on}~M_-\bigr\}.
\]
\end{thm}

\begin{rem}
Our proof of the general spectral flow formula bases on the property
{\bf (a1')} and {\bf (a2)}, and by making use
of the general spectral flow formula and $L_2$ reduction
theorem we prove the splitting formula above (\cite{FO}).
\end{rem}

Let $\mathfrak{D}_{APS}^{\,\,0}$ be the space defined
in (\ref{aps-0}) for $A$ replaced by $A+C_0$ and denote
by $\mathfrak{D}_{APS}^{\,\,1}$ the similar space
\begin{align}\label{aps-1}
&\mathfrak{D}_{APS}^{\,\,1}\\
&=\bigr\{ 
f\in H^1(M_+,\mathbb{E}_{|M_+})~\bigr|
~\text{if}~f_{\,|\Sigma}
=\sum\limits_{k \in\mathbb{Z}\backslash\{0\}}
a_k\psi_k,~\text{then}~a_{k}=0 ~\text{for}~k=1,2,3,\cdots
 \bigr\}.\notag
\end{align}
Note that the sections $\{\psi_k\}$ are now orthonormal 
eigensections of the tangential operator $B_1$ in the product form
\begin{equation}
A+C_1=\sigma\Bigr(\frac{\partial}{\partial u}+B_1\Bigr)
\end{equation}
corresponding to the parameter $t=1$ and 
should be chosen in such a way as 
noted in Remark (\ref{zero-mode}).

We have  continuous curves 
$\bigr\{\gamma^{-}(\Ker (\mathcal{A}_{-}^{\,\,\,*}+C_t))\bigr\}$ 
of Cauchy data spaces in the Fredholm-Lagrangian-Grassmannian 
$\mathcal{F}\Lambda_{\mathfrak{D}_0}(\beta^{-})$ 
= $\mathcal{F}\Lambda_{\mathfrak{D}_{APS}^{\,\,0}}(\beta^{-})$ and\\ 
$\bigr\{\gamma^{+}(\Ker (\mathcal{A}_{+}^{\,\,\,*}+C_t))\bigr\}$
in $\mathcal{F}\Lambda_{\mathfrak{D}_1}(\beta^{+})$ 
= $\mathcal{F}\Lambda_{\mathfrak{D}_{APS}^{\,\,1}}(\beta^{+})$. 

The H\"ormander index
is defined for four Lagrangian subspaces
$\mu_0$, 
$\mu_1$,
$\gamma^-(\mathfrak{D}_0)$ 
and $\gamma^-(\mathfrak{D}_{APS}^{\,\,0})$,
where
$\mu_i\in\mathcal{F}\Lambda_{\mathfrak{D}_0}(\beta^-)
=\mathcal{F}\Lambda_{\mathfrak{D}_{APS}^{\,\,0}}(\beta^-)$,
also defined for $\nu_0$, $\nu_1$,
$\gamma^+(\mathfrak{D}_1)$ 
and $\gamma^+(\mathfrak{D}_{APS}^{\,\,1})$
$\bigr(\nu_i\in\mathcal{F}\Lambda_{\mathfrak{D}_1}(\beta^+)=
\mathcal{F}\Lambda_{\mathfrak{D}_{APS}^{\,\,1}}(\beta^+)$$\bigr)$,
as noted in Remark (\ref{flg}).  

Since 
\begin{align*}
&\Sf\bigr(\{A_{\mathfrak{D}_0}+C_t\}\bigr)\\
&=\Mas\bigr(\{\gamma^-(\Ker(\mathcal{A}_{-}^{\,\,\,*}+C_t))\},
\gamma^-(\mathfrak{D}_0)\bigr)\\
&=\Mas\bigr(\{\gamma^-(\Ker(\mathcal{A}_{-}^{\,\,\,*}+C_t))\},
\gamma^-(\mathfrak{D}_{APS}^{\,\,0})\bigr)\\
&\qquad\qquad
+\sigma_H\bigr(\gamma^-(\Ker(\mathcal{A}_-^{\,\,\,*}+C_0)),
\gamma^-(\Ker(\mathcal{A}_-^{\,\,\,*}+C_1));
\gamma^-(\mathfrak{D}_0),\gamma^-(\mathfrak{D}_{APS}^{\,\,0})\bigr)
\end{align*}
we have
\begin{thm}
\begin{align*}
&\Sf\bigr(\{A+C_t\}\bigr)\\
&=\Mas\bigr(\{\gamma^-(\Ker(\mathcal{A}_{-}^{\,\,\,*}+C_t))\},
\gamma^-(\mathfrak{D}_{APS}^{\,\,0})\bigr)\\
&\qquad\qquad 
+\sigma_H\bigr(\gamma^-(\Ker(\mathcal{A}_-^{\,\,\,*}+C_0)),
\gamma^-(\Ker(\mathcal{A}_-^{\,\,\,*}+C_1))
;\gamma^-(\mathfrak{D}_0),\gamma^-(\mathfrak{D}_{APS}^{\,\,0})\bigr)\\
&+\Mas\bigr(\{\gamma^+(\Ker({\mathcal{A}_{+}}^{\,\,\,*}+C_t))\},
\gamma^+(\mathfrak{D}_{APS}^{\,\,1})\bigr)\\
&\qquad\qquad
+\sigma_H\bigr(\gamma^+(\Ker(\mathcal{A}_+^{\,\,\,*}+C_0)),
\gamma^+(\Ker(\mathcal{A}_+^{\,\,\,*}+C_1))
;\gamma^+(\mathfrak{D}_1),
\gamma^+(\mathfrak{D}_{APS}^{\,\,1})\bigr)\\
&=\Sf\bigr(\{A_{\mathfrak{D}_{APS}^{\,\,0}}+C_t\}\bigr) +
\Sf\bigr(\{A_{\mathfrak{D}_{APS}^{\,\,1}}+C_t\}\bigr)\\
&\qquad\qquad
+\sigma_H\bigr(\gamma^-(\Ker(\mathcal{A}_-^{\,\,\,*}+C_0)),
\gamma^-(\Ker(\mathcal{A}_-^{\,\,\,*}+C_1))
;\gamma^-(\mathfrak{D}_0),\gamma^-(\mathfrak{D}_{APS}^{\,\,0})\bigr)\\
&\qquad\qquad+
\sigma_H\bigr(\gamma^+(\Ker(\mathcal{A}_+^{\,\,\,*}+C_0)),
\gamma^+(\Ker(\mathcal{A}_+^{\,\,\,*}+C_1))
;\gamma^+(\mathfrak{D}_1),\gamma^+(\mathfrak{D}_{APS}^{\,\,1})\bigr).
\end{align*}
\end{thm}

\begin{cor}
If the family $\{A+C_t\}$ is a loop, i.e., $C_0=C_1$, then we have
\begin{align}
&\Sf\bigr(\{A+C_t\}\bigr)=\Sf\bigr(\{A_{\mathfrak{D}_{APS}^{\,\,0}} + C_t\}\bigr) 
+\Sf\bigr(\{A_{\mathfrak{D}_{APS}^{\,\,1}} + C_t\}\bigr).
\end{align}
\end{cor}

\begin{rem}
Although it holds the spectral flow formula expressed
in terms of the Maslov index
of Cauchy data spaces under the assumption {\bf (a1')}, 
it would not be clear whether the splitting formulas
of the spectral flow like above formulas hold always without the second
assumption {\bf (a2)}. 
Such assumptions are fit
to the framework of the symplectic Hilbert space theory,
after once the spaces $\beta^{\pm}$ are determined.
However it would be expected that generalizations of
splitting formula of spectral flow and the index similar to
(\ref{asymmetry}) without the assumption {\bf (a2)}
would be carried out through a further analysis
of the pseudo-differential operator
theory including the Calder\'on projector
and the operator $\mathcal{P}_{\tau(\Ker(\mathcal{A}_+^{\,\,\,*}))}$.
\end{rem}


\end{document}